# BERGMAN-TYPE PROJECTIONS IN GENERALIZED FOCK SPACES

Hélène Bommier-Hato, Miroslav Engliš and El-Hassan Youssfi

ABSTRACT. We give criteria for boundedness of the associated Bergman-type projections on $L^p$ spaces on $\mathbf{C}^n$ with respect to generalized Gaussian weights $\exp(-|z|^{2m})$, $m > 0$. Complete characterization is obtained for $0 < m \leq 1$ and $\frac{2n}{2n-1} < m < 2$; partial results are given for other values of $m$.

## 1. Introduction

Let $H_\alpha^p = H_\alpha^p(\mathbf{C}^n)$ be the space of all entire functions $f$ on $\mathbf{C}^n$, $n \geq 1$, such that $|f|^p$ is integrable with respect to the Gaussian

$$d\mu_\alpha(z) := e^{-\alpha|z|^2} \, dz, \tag{1}$$

where $dz$ stands for the Lebesgue volume on $\mathbf{C}^n$ and $\alpha > 0$, $1 \leq p < \infty$. Equipped with the norm inherited from $L_\alpha^p = L_\alpha^p(\mathbf{C}^n, d\mu_\alpha)$, $H_\alpha^p$ become Banach spaces; in particular, $H_\alpha^2$ is the Segal-Bargmann-Fock space of quantum mechanics with parameter $\alpha$ [11]. The function

$$K_\alpha(x,y) := \left(\frac{\alpha}{\pi}\right)^n e^{\alpha\langle x,y\rangle}, \qquad x,y \in \mathbf{C}^n,$$

is the reproducing kernel for $H_\alpha^2$, and the integral operator

$$P_\alpha f(x) := \int_{\mathbf{C}^n} f(y) \, K_\alpha(x,y) \, d\mu_\alpha(y), \qquad x \in \mathbf{C}^n,$$

is the orthogonal projection in $L_\alpha^2$ onto $H_\alpha^2$ (an analogue of the Bergman projection). By construction, $P_\alpha$ is bounded on $L_\alpha^2$, but this turns out to be no longer the case for $L_\alpha^p$ with $p \neq 2$. In fact, Janson, Peetre and Rochberg proved the following assertions.

---

1991 *Mathematics Subject Classification*. Primary 32A36; Secondary 47G10, 46E22, 33E12.

*Key words and phrases*. Bergman projection, higher order Fock space, reproducing kernel, Mittag-Leffler function.

Research supported by GA AV ČR grant no. IAA100190802, Czech Ministry of Education research plan no. MSM4781305904, and the French ANR DYNOP, Blanc07-198398.

**Theorem.** [16, §9.1] *Let $\alpha \in \mathbf{R}$, $\beta > 0$, $1 \leq p < \infty$ satisfy $\beta p > \alpha$. Then $P_\beta$ is bounded from $L^p_\alpha$ into $L^p_\gamma$, where $\frac{1}{\gamma} = \frac{4(\beta p - \alpha)}{p^2 \beta^2}$. In particular, $P_\beta$ is bounded on $L^p_\alpha$ if and only if $p\beta = 2\alpha$. Finally, $P_\beta$ maps $L^p_\beta$ into $L^q_\beta$ ($1 \leq q < \infty$) if and only if either $q < 4(1 - \frac{1}{p})$ or $p = q = 2$.*

The second part of the theorem was later proved by a different method by Dostanić and Zhu [6], who also computed the exact value of the norm of $P_\beta : L^p_\alpha \to L^p_\alpha$ with $2\alpha = p\beta$.

In this paper, we extend the last theorem in two ways: first of all, we discuss also the boundedness of $P_\beta$ from $L^p_\alpha$ into $L^q_\gamma$ for $p \neq q$ and $\alpha \neq \gamma$; and second, we consider the more general Fock-type spaces $H^p_{\alpha m} = H^p_{\alpha m}(\mathbf{C}^n)$ (and the associated projections $P_{\alpha m}$) consisting of all entire functions in $L^p_{\alpha m} := L^p(\mathbf{C}^n, d\mu_{\alpha m})$, where $d\mu_{\alpha m}$ are the generalized Gaussians

$$d\mu_{\alpha m}(z) := e^{-\alpha |z|^{2m}} dz, \qquad \alpha, m > 0,$$

which reduce to (1) for $m = 1$. We obtain a complete result for $m \leq 1$ and $\frac{2n}{2n-1} < m < 2$. Namely, denote

$$q_{\max} = q_{\max}(\alpha, \beta, \gamma, p) := \frac{4\gamma}{\beta^2}\Big(\beta - \frac{\alpha}{p}\Big).$$

**Theorem 1.** *Let $\alpha \in \mathbf{R}$, $\beta, \gamma > 0$, $1 \leq p, q < \infty$ and $0 < m \leq 1$. Then the following are equivalent:*
  (i) *$P_{\beta m} : L^p_{\alpha m} \to L^q_{\gamma m}$ is bounded;*
  (ii) *either $q < q_{\max}$, or $q = q_{\max} \geq p$.*

**Theorem 2.** *Let $\alpha \in \mathbf{R}$, $\beta, \gamma > 0$, $1 \leq p, q < \infty$ and $\frac{2n}{2n-1} < m < 2$. Then the following are equivalent:*
  (i) *$P_{\beta m} : L^p_{\alpha m} \to L^q_{\gamma m}$ is bounded;*
  (ii) *either $q < q_{\max}$, or $q = q_{\max} = p$.*

The cases of $1 < m \leq \frac{2n}{2n-1}$ and $m \geq 2$ seem to be more subtle, and the information we can offer is less complete. We prove the following:
  – If $q < q_{\max}$, then $P_{\beta m} : L^p_{\alpha m} \to L^q_{\gamma m}$ is bounded.
  – If $q > q_{\max}$, then $P_{\beta m} : L^p_{\alpha m} \to L^q_{\gamma m}$ is unbounded.
  – If $q = q_{\max} < p$ and $m < 2$, then $P_{\beta m} : L^p_{\alpha m} \to L^q_{\gamma m}$ is unbounded.
  – If $q = q_{\max} > p$ and $m > \frac{2n}{2n-1}$, then $P_{\beta m} : L^p_{\alpha m} \to L^q_{\gamma m}$ is unbounded.
  – If $q = q_{\max} = p$, then $P_{\beta m} : L^p_{\alpha m} \to L^q_{\gamma m}$ is bounded.

The general picture is thus a bit unclear at the moment, and we have no conjectures for the two remaining cases $q = q_{\max} > p$, $1 < m \leq \frac{2n}{2n-1}$, and $q = q_{\max} < p$, $m \geq 2$.

Note that for $\alpha = \gamma$, it is easy to see that always $q_{\max} \leq p$, with equality only for $p = 2\alpha/\beta$; thus our Theorem 1 indeed recovers the results of Dostanić and Zhu, i.e. the second part of the theorem by Janson, Peetre and Rochberg mentioned above. Similarly for the first and the third part of the theorem, i.e. for $p = q$, $\frac{1}{\gamma} = \frac{4(\beta p - \alpha)}{p^2 \beta^2}$ and $\alpha = \beta = \gamma$, respectively. (The case $m = 1$, $\alpha = \gamma$ of Theorem 1, i.e. the variant for $p \neq q$ of the result of Dostanić and Zhu, appears also in §5 of [7].)

We remark that $\alpha$ is allowed to be negative in the last two theorems; however, since $L^q_{\gamma m}$ contains no holomorphic functions for $\gamma \leq 0$, there is no point in considering nonpositive $\beta$ or $\gamma$.

The proof in [6] was based on the Schur test, which however turns out not to work very well for us here (even in its strengthened version due to Gagliardo [12]; see the remarks in Section 8.2 below.) Instead, we proceed by an interpolation argument, adapted from the approach in [16]. The crucial ingredient are estimates for various integrals of the corresponding reproducing kernels $K_{\alpha m}(x, y)$ of $H^2_{\alpha m}$, which involve interesting special functions (the generalized Mittag-Leffler function and its derivatives). Simpler versions of these estimates can be found in the papers [4], [5] by the first and the third author.

As noted above, Dostanić and Zhu in [6] also computed the exact norm of $P_\beta : L^p_\alpha \to L^p_\alpha$ in the case when it is bounded ($p = \frac{2\alpha}{\beta}$). We have made no effort to establish such results also in our more general setting, though conceivably this could be done in an analogous manner as in [6] at least in some cases.

The relevant prerequisites on the spaces $H^p_{\alpha m}$ and their reproducing kernels $K_{\alpha m}$ are assembled in Section 2, together with some handy formulas for norms of monomials. The latter are used in Section 4 to get a necessary condition for boundedness of $P_{\beta m}$ from $L^p_{\alpha m}$ to $L^q_{\gamma m}$, after making a convenient reduction of the problem in Section 3. A sufficient condition for boundedness is established in Section 5. The next two sections deal with the difficult case $q = q_{\max}$: miscellaneous results for general $m$ are obtained in Section 6, and Theorems 1 and 2 are proved in Section 7. The final Section 8 contains some concluding remarks and comments.

NOTATION. Throughout the paper, $\alpha, \beta, \gamma, p$ and $q$ are real numbers with $\beta, \gamma > 0$ and $1 \leq p, q < \infty$, and $\|\cdot\|_{p\alpha m}$ denotes the norm in $L^p_{\alpha m}$. For $p = 2$, we abbreviate $\|\cdot\|_{p\alpha m}$ to $\|\cdot\|_{\alpha m}$. The norm in $\mathbf{C}^n$ is denoted simply by $|\cdot|$ (even for $n > 1$). The notation $f = O(g)$, or $|f| \lesssim |g|$, for two functions $f, g$ on $\mathbf{C}^n$ means that there exists a constant $C$ such that $|f(z)| \leq C|g(z)|$ for all $z$ with $|z|$ large enough. If $f = O(g)$ and $g = O(f)$, we will write $f \sim g$.

ACKNOWLEDGEMENT. Large part of the work on this paper was done while the second author was visiting the other two; it is a great pleasure for him to gratefully acknowledge the hospitality and support of the LATP of Université de Provence in Marseille extended to him on that occasion.

## 2. Fock-type spaces

Since the weight function $e^{-\alpha |z|^{2m}}$ depends only on $|z|$, it follows in the standard way (integrating in polar coordinates) that the monomials $z^\nu$, $\nu$ a multiindex, form an orthogonal basis in $H^2_{\alpha m}$. Using the formula

$$(2) \qquad \int_{\mathbf{S}^{2n-1}} |\zeta^\nu|^2 \, d\sigma(\zeta) = \frac{2\nu! \pi^n}{\Gamma(n + |\nu|)}$$

for integration over the unit sphere $\mathbf{S}^{2n-1}$ in $\mathbf{C}^n$ with respect to the surface measure $d\sigma$, we obtain the expression for their norm squares

$$(3) \qquad \|z^\nu\|^2_{\alpha m} = \int_{\mathbf{C}^n} |z^\nu|^2 \, d\mu_{\alpha m}(z) = \frac{\pi^n \nu!}{m\Gamma(n + |\nu|)} \frac{\Gamma(\frac{n+|\nu|}{m})}{\alpha^{\frac{n+|\nu|}{m}}}.$$

Here we started using the usual multiindex notations $\nu! = \nu_1! \ldots \nu_n!$, $z^\nu = z_1^{\nu_1} \ldots z_n^{\nu_n}$ and so on. By the standard formula for a reproducing kernel [1], it follows that

$$K_{\alpha m}(x,y) = \sum_\nu \frac{x^\nu \overline{y}^\nu}{\|z^\nu\|_{\alpha m}^2} \tag{4}$$

$$= \frac{m\alpha^{n/m}}{\pi^n} \sum_{k=0}^\infty \frac{(\alpha^{1/m}\langle x,y\rangle)^k}{k!} \frac{\Gamma(n+k)}{\Gamma(\frac{n+k}{m})} \tag{5}$$

(see Lemma 5.1 in [5]; the formula there differs from ours, due to different normalizations, by a factor of $\frac{2\pi^n}{(n-1)!}$, which is immaterial for our purposes).

The corresponding orthogonal projection $P_{\alpha m} : L^2_{\alpha m} \to H^2_{\alpha m}$ is an integral operator with kernel $K_{\alpha m}$ (with respect to the measure $d\mu_{\alpha m}$):

$$P_{\alpha m} f(x) = \int_{\mathbf{C}^n} f(y) \, K_{\alpha m}(x,y) \, d\mu_{\alpha m}(y). \tag{6}$$

If $P_{\beta m}$ extends by continuity from $L^2_{\beta m} \cap L^p_{\alpha m}$ to all of $L^p_{\alpha m}$, then it is not a priori clear whether the extension will still be given by the formula (6). The next lemma shows that this is, however, the case at least in some special situations.

**Lemma 3.** *Assume that $P_{\beta m}$, defined by (6) on $L^2_{\beta m}$, extends to a continuous operator $L^p_{\alpha m} \to L^q_{\gamma m}$. For a multiindex $\nu$, let $L^{p,\nu}_{\alpha m}$ denote the subspace in $L^p_{\alpha m}$ of all functions of the form*

$$f(z) = z^\nu \phi(|z|). \tag{7}$$

*Then $|z^\nu|^2 \phi(z) \in L^1_{\beta m}$ and*

$$P_{\beta m} f(w) = \frac{w^\nu}{\|z^\nu\|_{\beta m}^2} \int_{\mathbf{C}^n} |z^\nu|^2 \, \phi(|z|) \, d\mu_{\beta m}(z). \tag{8}$$

*Proof.* The subset $L^{p,\nu}_{\alpha m} \cap L^{2,\nu}_{\beta m}$ is dense in $L^{p,\nu}_{\alpha m}$: it contains the functions $\iota_R(\chi_R f)$ for any $R > 0$, where $\chi_R$ is the characteristic function of the ball $|z| \leq R$, and $\iota_R(x) = x$ or $0$ according as $|x| \leq R$ or $|x| > R$; and $\iota_R(\chi_R f) \to f$ in $L^p_{\alpha m}$ as $R \to +\infty$. For $f \in L^{2,\nu}_{\beta m}$ we have by (6) and (4) (note that (4) converges in the $L^2_{\alpha m}$-norm for each fixed $y \in \mathbf{C}^n$) that $P_{\beta m} f(w)$ is given by the formula (8), where the integral exists since both $f$ and $z^\nu$ belong to $L^{2,\nu}_{\beta m}$.

Now if $f(z) = z^\nu \phi(|z|) \in L^{p,\nu}_{\alpha m}$, then also $g(z) := z^\nu |\phi(|z|)| \in L^{p,\nu}_{\alpha m}$; approximating the latter by functions $g_k$ in $L^{p,\nu}_{\alpha m} \cap L^{2,\nu}_{\beta m}$ as above, the continuity of $P_{\beta m} : L^p_{\alpha m} \to L^q_{\gamma m}$ implies that $P_{\beta m} g_k(w)$ has a finite limit as $k \to \infty$, for each $w$; by the Lebesgue Monotone Convergence Theorem, this means that $|z^\nu|^2 \phi(|z|) \in L^1_{\beta m}$, and approximating similarly $f$ by functions $f_k \in L^{p,\nu}_{\alpha m} \cap L^{2,\nu}_{\beta m}$ it follows by the Lebesgue Dominated Convergence Theorem that the validity of (8) for the $f_k$ implies that it is valid also for $f$. □

We will need a formula for $\|z^\nu\|_{p\alpha m}$ also for $p \neq 2$; though this is probably well-known, it is difficult to pinpoint a proof in the literature, so we include it here.

**Lemma 4.** *For any $\nu_1, \dots, \nu_n \geq 0$,*

$$\int_{\mathbf{S}^{2n-1}} |\zeta_1|^{2\nu_1} \dots |\zeta_n|^{2\nu_n} \, d\sigma(\zeta) = \frac{2\pi^n \nu!}{\Gamma(n+|\nu|)},$$

*where $|\nu| = \nu_1 + \dots + \nu_n$ and $\nu! = \Gamma(\nu_1+1)\dots\Gamma(\nu_n+1)$.*

Of course, for integer $\nu_1, \dots, \nu_n$ the last formula reduces to (2).

*Proof.* Let us denote the integral by $C_{2\nu}$. We express

$$(9) \qquad C_{2\nu} \int_0^\infty r^{2|\nu|} e^{-r^2} r^{2n-1} \, dr$$

in two ways. On the one hand, changing the variable to $t = r^2$ shows that this is equal to

$$C_{2\nu} \int_0^\infty t^{|\nu|+n-1} e^{-t} \frac{dt}{2} = \frac{1}{2}\Gamma(n+|\nu|) C_{2\nu}.$$

On the other hand, (9) is obtained upon making the change of variable $z = r\zeta$ in the integral

$$\int_{\mathbf{C}^n} |z_1|^{2\nu_1} \dots |z_n|^{2\nu_n} e^{-|z|^2} \, dz,$$

which equals

$$\prod_{j=1}^n \int_{\mathbf{C}} |z_j|^{2\nu_j} e^{-|z_j|^2} \, dz_j = \prod_{j=1}^n \int_0^{2\pi} d\theta \int_0^\infty r^{2\nu_j} e^{-r^2} r \, dr = \pi^n \nu!.$$

Comparing both formulas, the assertion follows. $\square$

**Corollary 5.** *For any multiindex $\nu$ and $p > 0$,*

$$(10) \qquad \|z^\nu\|_{p\alpha m}^p = \frac{\pi^n}{m} \frac{(\frac{\nu p}{2})!}{\Gamma(n+\frac{|\nu|p}{2})} \frac{\Gamma(\frac{|\nu|p+2n}{2m})}{\alpha^{\frac{|\nu|p+2n}{2m}}}.$$

*Proof.* Immediate from the previous lemma and integration in polar coordinates. $\square$

## 3. A REDUCTION

Given $p, q$, let us denote

$$(11) \qquad c \equiv c_{p,q}(\alpha, \beta, \gamma) = \frac{4\gamma}{\beta^2 q}\left(\beta - \frac{\alpha}{p}\right).$$

In other words, $c = q_{\max}(\alpha, \beta, \gamma, p)/q$.

**Proposition 6.** *For fixed $p$ and $q$, the boundedness or unboundedness of $P_{\beta m}$ from $L^p_{\alpha m}$ into $L^q_{\gamma m}$ depends only on the value of $c$.*

*Proof.* For $\delta > 0$, consider the dilation operator

$$D_\delta f(z) := f(\delta^{1/m} z).$$

A simple change-of-variable computation shows that $D_\delta$ maps $L^p_{\alpha m}$ isomorphically onto $L^p_{\delta^2\alpha,m}$; in fact, $\delta^{2n/m}D_\delta$ is an isometry of the former onto the latter. Since, in view of (5), $K_{\beta m}(z,\delta^{-1/m}w) = \delta^{2n/m}K_{\beta/\delta^2,m}(\delta^{1/m}z,w)$, another simple change of variable argument shows that

$$P_{\beta m}D_\delta = D_\delta P_{\beta/\delta^2,m}.$$

Thus $P_{\beta m}$ is bounded $L^p_{\alpha m} \to L^q_{\gamma m}$ if and only if $P_{\beta/\delta^2,m}$ is bounded from $D_{1/\delta}L^p_{\alpha m} = L^p_{\alpha/\delta^2,m}$ into $D_{1/\delta}L^q_{\gamma m} = L^q_{\gamma/\delta^2,m}$. Let us temporarily call, for $p,q$ fixed, $(\alpha,\beta,\gamma)$ a *bounded* triple if $P_{\beta m} : L^p_{\alpha m} \to L^q_{\gamma m}$ is bounded. Then we have just shown that for any $1/\delta^2 = t > 0$,

(12) $\qquad (\alpha,\beta,\gamma)$ is bounded $\iff (t\alpha, t\beta, t\gamma)$ is bounded.

Next, for $\epsilon \in \mathbf{R}$ consider the operator

$$E_{\epsilon m}f(z) := e^{\epsilon|z|^{2m}}f(z).$$

Then $E_{\epsilon m}$ is an isometric isomorphism of $L^p_{\alpha m}$ onto $L^p_{\alpha+p\epsilon,m}$. Furthermore, using again (5),

$$P_{\delta m}E_{\delta-\beta,m}f(z) = \int_{\mathbf{C}^n} f(w)\, e^{(\delta-\beta)|w|^{2m}} K_{\delta m}(z,w)\, e^{-\delta|w|^{2m}}\, dw$$
$$= \int_{\mathbf{C}^n} f(w)K_{\delta m}(z,w)\, e^{-\beta|w|^{2m}}\, dw$$
$$= \int_{\mathbf{C}^n} f(w)\Big(\frac{\delta}{\beta}\Big)^{n/m} K_{\beta m}\Big(\frac{\delta^{1/m}}{\beta^{1/m}}z,w\Big)\, e^{-\beta|w|^{2m}}\, dw$$
$$= (\tfrac{\delta}{\beta})^{n/m} D_{\delta/\beta}P_{\beta m}f(z),$$

or

$$P_{\delta m}E_{\delta-\beta,m} = (\tfrac{\delta}{\beta})^{n/m}D_{\delta/\beta}P_{\beta m}$$

for any $\delta,\beta > 0$. It follows that $P_{\beta m} : L^p_{\alpha m} \to L^q_{\gamma m}$ is bounded if and only if $P_{\delta m}$ is bounded from $E_{\delta-\beta,m}L^p_{\alpha m} = L^p_{\alpha+(\delta-\beta)p,m}$ into $D_{\delta/\beta}L^q_{\gamma m} = L^q_{\delta^2\gamma/\beta^2,m}$. In other words, for any $\delta > 0$,

$$(\alpha,\beta,\gamma) \text{ is bounded} \iff (\alpha+(\delta-\beta)p, \delta, \frac{\delta^2\gamma}{\beta^2}) \text{ is bounded}.$$

Taking in particular $\delta = \beta^2 q/\gamma$, we get

$$(\alpha,\beta,\gamma) \text{ is bounded} \iff \Big(\alpha+p\Big(\frac{\beta^2 q}{\gamma}-\beta\Big), \frac{\beta^2 q}{\gamma}, \frac{\beta^2 q^2}{\gamma}\Big) \text{ is bounded}$$
$$\iff \Big(\frac{\alpha\gamma}{\beta^2 q}+p-\frac{p\gamma}{\beta q}, 1, q\Big) \text{ is bounded}$$

by (12). Since

$$\frac{\alpha\gamma}{\beta^2 q} - \frac{p\gamma}{\beta q} = \frac{\gamma p}{\beta^2 q}\Big(\frac{\alpha}{p}-\beta\Big) = -\frac{p}{4}c,$$

we get that

$$(\alpha, \beta, \gamma) \text{ is bounded} \iff ((1 - \tfrac{c}{4})p, 1, q) \text{ is bounded},$$

and the claim follows. $\square$

Since evidently $L^q_{\gamma m} \hookrightarrow L^q_{\delta m}$ if $\delta > \gamma$ and $L^p_{\epsilon m} \hookrightarrow L^p_{\alpha m}$ if $\epsilon < \alpha$, we have that if $(\alpha, \beta, \gamma)$ is bounded then so is $(\epsilon, \beta, \delta)$ for all $\epsilon < \alpha$ and $\delta > \gamma$. It follows that there exists $c_0 = c_0(p, q) \in \{-\infty\} \cup [0, +\infty)$ such that

$$(\alpha, \beta, \gamma) \text{ is bounded if } c > c_0 \text{ and unbounded if } c < c_0.$$

Similarly, since $\mu_{\gamma m}$ is a finite measure (as $\gamma > 0$), we have $L^q_{\gamma m} \hookrightarrow L^r_{\gamma m}$ if $r < q$ (and likewise $L^s_{\alpha m} \hookrightarrow L^p_{\alpha m}$ for $s > p$ if $\alpha > 0$); it transpires that $c_0(p, q)$ must be a nondecreasing function of $q$ and nonincreasing function of $p$. Finally, due to a simple consequence of the Hölder inequality to the effect that[1]

$$L^s_{\delta m} \hookrightarrow L^p_{\alpha m} \text{ if } s \geq p \text{ and } \frac{\alpha}{p} > \frac{\delta}{s},$$

one could draw further conclusions about $c_0(p, q)$ by letting all of $p, q, \alpha, \gamma$ vary. We will not pursue this line further, however, since we are now going to show directly that $c_0(p, q) = 1$ for all $p$ and $q$.

## 4. A NECESSARY CONDITION

We keep the notation (11) throughout the rest of this paper.

**Proposition 7.** *If $P_{\beta m} : L^p_{\alpha m} \to L^q_{\gamma m}$ is bounded, then $c \geq 1$.*

*Proof.* Consider the function

$$f(z) = z^\nu e^{-\lambda |z|^{2m}}$$

where $\nu$ is a multiindex and $\lambda < \alpha/p$. The condition on $\lambda$ guarantees that $f \in L^p_{\alpha m}$, with

$$\|f\|^p_{p\alpha m} = \int_{\mathbf{C}^n} |z^\nu|^p e^{-(\alpha - \lambda p)|z|^{2m}} dz$$
$$= \|z^\nu\|^p_{p, \alpha - \lambda p, m}.$$

---

[1] The assertion follows from Hölder's inequality applied to the pair of functions $|f(z)|^p e^{-b|z|^{2m}}$ and 1 with respect to the measure $e^{-c|z|^{2m}} dz$, $c > 0$, and exponent $\frac{s}{p}$:

$$\int |f(z)|^p e^{-b|z|^{2m}} e^{-c|z|^{2m}} dz \leq \left( \int |f(z)|^s e^{-bs|z|^{2m}/p} e^{-c|z|^{2m}} dz \right)^{p/s} \left( \int e^{-c|z|^{2m}} dz \right)^{1-(p/s)},$$

or

$$\|f\|_{p, b+c, m} \leq C \|f\|_{s, c+bs/p, m}, \qquad C = \left( \int e^{-c|z|^{2m}} dz \right)^{1/p - 1/s} < \infty.$$

The condition $\frac{\alpha}{p} > \frac{\delta}{s}$ guarantees that we can choose $c > 0$ and $b \in \mathbf{R}$ so that $b+c = \alpha$, $c+bs/p = \delta$. Hence $\|f\|_{p\alpha m} \leq C \|f\|_{s\delta m}$.

If $P_{\beta m} : L^p_{\alpha m} \to L^q_{\gamma m}$ is bounded, then by Lemma 3, $|z^\nu|^2 e^{-\lambda |z|^{2m}} \in L^1_{\beta m}$ and $P_{\beta m} f$ is given by (8). The former implies that $\lambda < \beta$, or, since $\lambda$ can be taken arbitrary $< \alpha/p$,

$$\tag{13} \frac{\alpha}{p} \leq \beta;$$

the latter then implies that

$$P_{\beta m} f(w) = \frac{w^\nu}{\|z^\nu\|^2_{\beta m}} \|z^\nu\|^2_{\beta-\lambda,m} = \left(\frac{\beta}{\beta-\lambda}\right)^{\frac{n+|\nu|}{m}} w^\nu$$

by (3).

Altogether, we see that if $P_{\beta m}$ extends to a bounded operator from $L^p_{\alpha m}$ into $L^q_{\gamma m}$ — with norm $C$, say — then $\frac{\alpha}{p} \leq \beta$ and

$$\tag{14} \left(\frac{\beta}{\beta-\lambda}\right)^{\frac{n+|\nu|}{m}} \|z^\nu\|_{q\gamma m} = \|P_{\beta m} f\|_{q\gamma m} \leq C \|f\|_{p\alpha m} = C \|z^\nu\|_{p,\alpha-\lambda p, m},$$

for all multiindices $\nu$ and $\lambda < \frac{\alpha}{p}$.

We now specialize to $\nu = (k, 0, \ldots, 0)$, $k = 0, 1, 2, \ldots$, so that

$$\|z^\nu\|_{p\alpha m} = \left[\frac{\pi^n}{m} \frac{\Gamma(\frac{kp}{2}+1)}{\Gamma(\frac{kp}{2}+n)} \frac{\Gamma(\frac{kp+2n}{2m})}{\alpha^{\frac{kp+2n}{2m}}}\right]^{1/p}$$

by (10). Then (14) becomes

$$\left(\frac{\beta}{\beta-\lambda}\right)^{\frac{k+n}{m}} \left[\frac{\pi^n}{m} \frac{\Gamma(\frac{kq}{2}+1)}{\Gamma(\frac{kq}{2}+n)} \frac{\Gamma(\frac{kq+2n}{2m})}{\gamma^{\frac{kq+2n}{2m}}}\right]^{1/q} \leq C \left[\frac{\pi^n}{m} \frac{\Gamma(\frac{kp}{2}+1)}{\Gamma(\frac{kp}{2}+n)} \frac{\Gamma(\frac{kp+2n}{2m})}{(\alpha-\lambda p)^{\frac{kp+2n}{2m}}}\right]^{1/p}$$

for all $k = 0, 1, 2, \ldots$ and $\lambda < \frac{\alpha}{p}$, with some $C$ independent of $\lambda$ and $k$. Now raise both sides to the power $\frac{2m}{k}$, divide by $k$, and let $k \to \infty$. Since, for any $\rho > 0$ and $\sigma \in \mathbf{R}$,

$$\lim_{k\to\infty} \frac{\Gamma(\rho k + \sigma)^{1/\rho k}}{k} = \frac{\rho}{e}$$

(by Stirling's formula), we obtain

$$\left(\frac{\beta}{\beta-\lambda}\right)^2 \frac{q/2me}{\gamma} \leq \frac{p/2me}{\alpha-\lambda p},$$

for all $\lambda < \frac{\alpha}{p}$. Cancelling $2me$ and setting $\lambda = \frac{\alpha}{p} - t$, this becomes

$$\frac{\beta^2}{(\beta - \frac{\alpha}{p} + t)^2} \frac{q}{\gamma} \leq \frac{1}{t}$$

or

$$\tag{15} -\frac{\beta^2}{\gamma} qt + \left(\beta - \frac{\alpha}{p} + t\right)^2 \geq 0 \qquad \forall t > 0.$$

The left-hand side attains its minimum at $t = t_{\min} := \frac{\beta^2}{2\gamma}q - \beta + \frac{\alpha}{p}$, and equals $(\beta - \frac{\alpha}{p})^2 \geq 0$ at $t = 0$. If $t_{\min} \leq 0$, the inequality in (15) thus indeed holds for all positive $t$. If $t_{\min} > 0$, then (15) holds if and only if the inequality there holds for $t = t_{\min}$, that is, if and only if

$$(\beta - \tfrac{\alpha}{p})^2 - t_{\min}^2 \geq 0$$

or (since $t_{\min} > 0$ and $\beta \geq \frac{\alpha}{p}$)

$$0 < t_{\min} \leq \beta - \tfrac{\alpha}{p}.$$

Altogether, we see that (15) is equivalent to $t_{\min} \leq \beta - \frac{\alpha}{p}$, or

$$q \leq \frac{4\gamma}{\beta^2}\left(\beta - \frac{\alpha}{p}\right) = q_{\max}(\alpha, \beta, \gamma, p),$$

or $c = q_{\max}/q \geq 1$, completing the proof. □

Note that since $q \geq 1$, the boundedness of $P_{\beta m} : L^p_{\alpha m} \to L^q_{\gamma m}$ implies that

(16) $$\beta - \frac{\alpha}{p} > 0,$$

i.e. we must even have strict inequality in (13).

## 5. A sufficient condition

In view of (5), the reproducing kernel $K_{\alpha m}$ can be expressed in terms of the Mittag-Leffler function

$$E_{\alpha,\beta}(z) := \sum_{k=0}^{\infty} \frac{z^k}{\Gamma(\alpha k + \beta)}, \qquad \alpha, \beta > 0,$$

as

(17) $$K_{\alpha m}(x, y) = \frac{m\alpha^{n/m}}{\pi^n} E^{(n-1)}_{\frac{1}{m},\frac{1}{m}}(\alpha^{1/m}\langle x, y\rangle).$$

From the formula (22) in vol. III, §18.1 of [3], one learns that $E_{1/m,1/m}(z)$ has the asymptotic expansion

$$E_{\frac{1}{m},\frac{1}{m}}(z) = mz^{m-1}e^{z^m}[1 + O(e^{-\theta|z|^m/\sqrt{2}})]$$

as $z \to \infty$, $|\arg z| \leq \frac{\pi}{4m}$, with $\theta > 0$. (One can take $\theta = (1 - \cos 2\pi m)$ for $0 < m < \frac{1}{4}$, and any $0 < \theta < 1$ for $m \geq \frac{1}{4}$.) Furthermore, the expansion can be differentiated termwise any number of times: this can be seen either by checking the derivation from the contour integral representation in [3], or by using Cauchy's formula (integrating over circles of radius 1 around $z$). In particular, we have

$$E^{(n-1)}_{\frac{1}{m},\frac{1}{m}}(z) = O(z^{(m-1)n}e^{z^m})$$

as $z \to +\infty$, $z > 0$. Since $|E^{(n-1)}_{\frac{1}{m},\frac{1}{m}}(z)| \leq E^{(n-1)}_{\frac{1}{m},\frac{1}{m}}(|z|)$ (the power series for $E^{(n-1)}_{\frac{1}{m},\frac{1}{m}}$ has nonnegative coefficients), we see that

(18) $$|K_{\alpha m}(x, y)| \lesssim |\langle x, y\rangle|^{(m-1)n}e^{\alpha|\langle x,y\rangle|^m}$$
$$\lesssim (|x|\,|y|)^{(m-1)n}e^{\alpha|x|^m|y|^m}$$

as $|x||y| \to +\infty$. (See [4] for a detailed proof.)

We will need a simple computational lemma.

**Lemma 8.** *For $k > -2n$, $B > 0$ and $A > 0$,*
$$\int_{\mathbf{C}^n} |w|^k \, e^{B|w|^m} e^{-A|w|^{2m}} \, dw \lesssim B^{\frac{k+2n}{m}-1} e^{B^2/4A}$$
*as $B \to +\infty$.*

*Proof.* Denote, for brevity, $\frac{k+2n}{m} - 1 =: \rho$; by hypothesis, $\rho > -1$. Passing to polar coordinates $w = r\zeta$ shows that
$$\int_{\mathbf{C}^n} |w|^k \, e^{B|w|^m} e^{-A|w|^{2m}} \, dw = \frac{2\pi^n}{\Gamma(n)} \int_0^\infty r^{k+2n-1} e^{Br^m - Ar^{2m}} \, dr$$
$$= \frac{2\pi^n}{\Gamma(n)m} \int_0^\infty t^\rho \, e^{Bt - At^2} \, dt.$$

The last integral can be evaluated explicitly in terms of the parabolic cylinder functions $D_\nu$, namely,
$$\int_0^\infty t^\rho \, e^{Bt - At^2} \, dt = \frac{\Gamma(\rho+1)}{(2A)^{(\rho+1)/2}} e^{B^2/8A} D_{-\rho-1}\Big(\frac{-B}{\sqrt{2A}}\Big).$$

See [3], formula (3) in vol. II, §8.3, or formula 3.462.1 in [14]. Using the known asymptotic behaviour of $D_\nu$ (see [20], §16.52; this is also reproduced in vol. II, §8.4 of [3], but with missing parentheses!)
$$D_\nu(z) \sim (-z)^{-\nu-1} e^{z^2/4} \qquad \text{as } z \to -\infty,$$
we obtain the lemma (even with $\lesssim$ replaced by $\sim$). □

*Remark.* It is not difficult to give a direct proof avoiding the parabolic cylinder functions; we indicate one for $\rho \geq 0$. (The case of $\rho < 0$ requires more labour.) We can continue the computation above by
$$\int_0^\infty t^\rho \, e^{Bt - At^2} \, dt = e^{B^2/4A} \int_{-B/2A}^\infty \Big(s + \frac{B}{2A}\Big)^\rho e^{-As^2} \, ds$$
$$= e^{B^2/4A} A^{-\rho-1} \int_{-B/2\sqrt{A}}^\infty \Big(t + \frac{B}{2\sqrt{A}}\Big)^\rho e^{-t^2} \, dt.$$

We need to show that the last integral is $O(B^\rho)$. Since $\rho \geq 0$,
$$\Big|t + \frac{B}{2\sqrt{A}}\Big|^\rho \leq \Big(|t| + \frac{B}{2\sqrt{A}}\Big)^\rho \leq C_\rho \Big[t^\rho + \Big(\frac{B}{2\sqrt{A}}\Big)^\rho\Big],$$
where $1/C_\rho$ denotes the (finite and positive) minimum of the function $x^\rho + (1-x)^\rho$ on the interval $[0,1]$. The last integral is therefore dominated by
$$C_\rho \int_{-\infty}^\infty \Big[|t|^\rho + \Big(\frac{B}{2\sqrt{A}}\Big)^\rho\Big] e^{-t^2} \, dt$$
$$= C_\rho\Big[\Gamma\Big(\frac{\rho+1}{2}\Big) + \Big(\frac{B}{2\sqrt{A}}\Big)^\rho \Gamma\Big(\frac{1}{2}\Big)\Big]$$
$$= O(B^\rho),$$
as asserted. □

The next proposition proves actually a little more than the boundedness of $P_{\beta m}$: it shows that $P_{\beta m}$ is, under the given hypothesis, still given by the formula (6) even for all $f \in L^p_{\alpha m}$ (i.e. not only for $f \in L^p_{\alpha m} \cap L^2_{\beta m}$).

**Proposition 9.** *If $c > 1$, then (6) defines a bounded operator from $L^p_{\alpha m}$ into $L^q_{\gamma m}$.*

*Proof.* By Hölder's inequality,

$$|P_{\beta m} f(z)| \leq \|f\|_{p\alpha m} \|K_{\beta m}(\,\cdot\,, z) e^{(\alpha-\beta)|\,\cdot\,|^{2m}}\|_{p'\alpha m}$$

where $\frac{1}{p} + \frac{1}{p'} = 1$. By (18),

$$\|K_{\beta m}(\,\cdot\,, z) e^{(\alpha-\beta)|\,\cdot\,|^{2m}}\|_{p'\alpha m}^{p'}$$
$$\lesssim \int_{\mathbf{C}^n} \left[ |z|^{(m-1)n} |w|^{(m-1)n} e^{\beta|z|^m |w|^m} e^{(\alpha-\beta)|w|^{2m}} \right]^{p'} e^{-\alpha|w|^{2m}} \, dw$$
$$= |z|^{(m-1)np'} \int_{\mathbf{C}^n} |w|^{(m-1)np'} \, e^{\beta p' |z|^m |w|^m - (\alpha - p'\alpha + p'\beta)|w|^{2m}} \, dw$$
$$\lesssim |z|^{2(m-1)np' + 2n - m} e^{\beta^2 p'^2 |z|^{2m} / 4(\alpha - p'\alpha + p'\beta)}$$

by the last lemma (note that $(m-1)np' \geq (m-1)n > -n > -2n$, and $\alpha - p'\alpha + p'\beta = p'(\beta - \frac{\alpha}{p}) = \frac{\beta^2 p'}{4\gamma} cq > 0$). As

$$\|P_{\beta m} f\|_{q\gamma m}^q \leq \|f\|_{p\alpha m}^q \int_{\mathbf{C}^n} \|K_{\beta m}(\,\cdot\,, z) e^{(\alpha-\beta)|\,\cdot\,|^{2m}}\|_{p'\alpha m}^q \, d\mu_{\gamma m}(z),$$

we see that (6) defines a bounded operator $L^p_{\alpha m} \to L^q_{\gamma m}$ if

$$(19) \quad \int_{\mathbf{C}^n} \left[ |z|^{2(m-1)np' + 2n - m} e^{\beta^2 p'^2 |z|^{2m} / 4(\alpha - p'\alpha + p'\beta)} \right]^{\frac{q}{p'}} e^{-\gamma|z|^{2m}} \, dz < \infty,$$

hence, if

$$\frac{\beta^2 p' q}{4(\alpha - p'\alpha + p'\beta)} - \gamma < 0,$$

or

$$q < \frac{4\gamma(\alpha - p'\alpha + p'\beta)}{\beta^2 p'} = \frac{4\gamma}{\beta^2}\left(\beta - \frac{\alpha}{p}\right) = q_{\max},$$

or $c = \frac{q_{\max}}{q} > 1$, completing the proof. $\square$

## 6. Further analysis

We are thus left with the borderline case $c = 1$. As announced in the Introduction, we have a complete characterizations only for $0 < m \leq 1$ and $\frac{2n}{2n-1} < m < 2$, which we will prove in the next section; here we present what we can say for the case of general $m$.

**Proposition 10.** *Assume that $P_{\beta m} : L^p_{\alpha m} \to L^q_{\gamma m}$ is bounded and $c = 1$. Then $(m-2)(p-q) \geq 0$.*

*Proof.* As in the proof of Proposition 7, consider the test functions

$$f(z) = z^\nu e^{\lambda |z|^{2m}},$$

where we now take $\lambda = 2\frac{\alpha}{p} - \beta$ and $\nu = (k, \ldots, k)$, $k = 0, 1, 2, \ldots$. Note that $\lambda < \frac{\alpha}{p}$ in view of (16). As in Section 4, the boundedness of $P_{\beta m} : L^p_{\alpha m} \to L^q_{\gamma m}$ implies that

$$\Big(\frac{\beta}{\beta - \lambda}\Big)^{\frac{|\nu|+n}{m}} \|z^\nu\|_{q\gamma m} \leq C \|z^\nu\|_{p, \alpha - \lambda p, m},$$

or

$$\Big(\frac{\beta}{\beta - \lambda}\Big)^{\frac{(k+1)n}{m}} \Big[\frac{\Gamma(\frac{kq}{2}+1)^n}{\Gamma(\frac{nkq}{2}+n)} \frac{\Gamma(\frac{nkq+2n}{2m})}{\gamma^{\frac{nkq+2n}{2m}}}\Big]^{\frac{1}{q}} \leq C \Big[\frac{\Gamma(\frac{kp}{2}+1)^n}{\Gamma(\frac{nkp}{2}+n)} \frac{\Gamma(\frac{nkp+2n}{2m})}{(\alpha - \lambda p)^{\frac{nkp+2n}{2m}}}\Big]^{\frac{1}{p}}$$

with some $C$ independent of $k$. Taking logarithms and noting that $\frac{\beta}{\beta - \lambda} = \frac{\beta p/2}{\beta p - \alpha}$ and $\alpha - \lambda p = \beta p - \alpha$, this yields

$$
\begin{aligned}
(20) \quad & \frac{(k+1)n}{m} \log \frac{\beta p/2}{\beta p - \alpha} + \frac{n}{2} \log \Gamma\Big(\frac{kq}{2}+1\Big)^{\frac{2}{q}} - \frac{n}{2} \log \Gamma\Big(\frac{nkq}{2}+n\Big)^{\frac{2}{nq}} \\
& + \frac{n}{2m} \log \Gamma\Big(\frac{nkq}{2m}+\frac{n}{m}\Big)^{\frac{2m}{nq}} - \frac{nkq+2n}{2mq} \log \gamma \\
& \leq C + \frac{n}{2} \log \Gamma\Big(\frac{kp}{2}+1\Big)^{\frac{2}{p}} - \frac{n}{2} \log \Gamma\Big(\frac{nkp}{2}+n\Big)^{\frac{2}{np}} \\
& + \frac{n}{2m} \log \Gamma\Big(\frac{nkp}{2m}+\frac{n}{m}\Big)^{\frac{2m}{np}} - \frac{nkp+2n}{2mp} \log(\beta p - \alpha).
\end{aligned}
$$

By Stirling's formula (cf. [3], vol. I, formula (3) in §1.18), we have for any $A > 0$ and $B \in \mathbf{R}$ as $k \to +\infty$

$$\log \Gamma(Ak+B)^{\frac{1}{A}} = \Big(k + \frac{B - \frac{1}{2}}{A}\Big)(\log A + \log k) - k + \frac{\log \sqrt{2\pi}}{A} + O\Big(\frac{1}{k}\Big).$$

Applying this to the logarithms of the Gamma functions in (20) and bringing everything except $C$ to one side of the inequality, the terms containing $k \log k$ cancel, while the coefficients at $k$ combine into

$$-\frac{n}{2m} \log \frac{4\gamma(\beta p - \alpha)}{\beta^2 pq} = -\frac{n}{2m} \log c,$$

which vanishes too since $c = 1$ by hypothesis. We thus end up with a relation of the form

$$(21) \qquad a \log k + b + O(\tfrac{1}{k}) \leq C,$$

with some constants $a, b$, for the first of which we get after some computations

$$a = \frac{(m-2)n(q-p)}{2mpq}.$$

Letting $k$ tend to infinity, (21) can only hold if $a \leq 0$, or

$$(m-2)(p-q) \geq 0,$$

proving our claim. $\square$

**Corollary 11.** *If $q = q_{\max} < p$ and $m < 2$, then $P_{\beta m}$ cannot be bounded from $L^p_{\alpha m}$ into $L^q_{\gamma m}$.*

The next proposition can be proved in exactly the same way as the last, only taking $\nu = (k, 0, \ldots, 0)$ (i.e. the same multiindex as in the proof of Proposition 7) and again $\lambda = 2\frac{\alpha}{p} - \beta$; we therefore omit the proof.

**Proposition 12.** *Assume that $P_{\beta m} : L^p_{\alpha m} \to L^q_{\gamma m}$ is bounded and $c = 1$. Then $((2n-1)m - 2n)(p-q) \geq 0$.*

**Corollary 13.** *If $q = q_{\max} > p$ and $m > \frac{2n}{2n-1}$, then $P_{\beta m}$ cannot be bounded from $L^p_{\alpha m}$ into $L^q_{\gamma m}$.*

Since $\frac{2n}{2n-1} \leq 2$, we also get

**Corollary 14.** *If $q = q_{\max} \neq p$ and $\frac{2n}{2n-1} < m < 2$, then $P_{\beta m}$ cannot be bounded from $L^p_{\alpha m}$ into $L^q_{\gamma m}$.*

Our wild guess for the borderline situation $c = 1$ (i.e. $q = q_{\max}$) is that $P_{\beta m} : L^p_{\alpha m} \to L^q_{\gamma m}$ is bounded if $0 < m \leq 1$ and $p \leq q$, or $1 < m < 2$ and $p = q$, or $m \geq 2$ and $p \geq q$, and unbounded in all other cases; however, we cannot procure too much further evidence towards such conclusion.

## 7. Proofs of Theorems 1 and 2

The results in this section are based on a sharper estimate for the integrals of reproducing kernels than can be obtained from (18) and Lemma 9.

Recall that the asymptotics of the Mittag-Leffler function $E_{1/m,1/m}(z)$ as $|z| \to +\infty$ are given by

$$E_{\frac{1}{m},\frac{1}{m}}(z) = \begin{cases} mz^{m-1}e^{z^m} + O(\frac{1}{z}), & |\arg z| \leq \frac{\pi}{2m}, \\ O(\frac{1}{z}), & \frac{\pi}{2m} < |\arg z| \leq \pi \end{cases}$$

for $m > \frac{1}{2}$, and by

$$E_{\frac{1}{m},\frac{1}{m}}(z) = m \sum_{j=-N}^{N} z^{m-1} e^{2\pi i j(m-1)} e^{z^m e^{2\pi i j m}} + O(\tfrac{1}{z}), \quad -\pi < \arg z \leq \pi,$$

for $0 < m \leq \frac{1}{2}$, where $N$ is the integer satisfying $N < \frac{1}{2m} \leq N+1$, and the powers $z^{m-1}, z^m$ are the principal branches. See e.g. Bateman and Erdelyi [3], vol. III, §18.1, formulas (21)–(22); other good references are Paris [17], Wong and Zhao [21], and §5.1.4 in the book by Paris and Kaminski [18], as well as Gorenflo, Loutchko and Luchko [13] and Hilfer and Seybold [15] which give some numerical data and illustrative graphics. (All these references also discuss in detail the Stokes phenomenon in the change of the asymptotics on the rays $|\arg z| = \frac{\pi}{m}$ for $m > 1$.)

As explained at the beginning of Section 4, the asymptotic expansions above can be differentiated termwise. Since a simple induction argument shows that

$$\frac{d^{k-1}}{dz^{k-1}}(mz^{m-1}e^{z^m}) = \frac{p_k(z^m)}{z^k}e^{z^m},$$

where $p_k$ are polynomials of degree $k$ defined recursively by

$$p_0 = 1, \qquad p_{k+1}(x) = (mx - k)p_k(x) + mxp'_k(x),$$

we thus obtain, in particular, the following asymptotics for $E^{(n-1)}_{1/m,1/m}(z)$ as $z \to \infty$:

$$E^{(n-1)}_{\frac{1}{m},\frac{1}{m}}(z) = \begin{cases} \dfrac{p_n(z^m)}{z^n} e^{z^m} + O(z^{-n}), & |\arg z| \leq \frac{\pi}{2m}, \\ O(z^{-n}), & \frac{\pi}{2m} < |\arg z| \leq \pi \end{cases}$$

for $m > \frac{1}{2}$, and

$$(22) \qquad E^{(n-1)}_{\frac{1}{m},\frac{1}{m}}(z) = \sum_{j=-N}^{N} \frac{p_n(z^m e^{2\pi ijm})}{z^n} e^{z^m e^{2\pi ijm}} + O(z^{-n}), \quad -\pi < \arg z \leq \pi,$$

for $0 < m \leq \frac{1}{2}$, $N < \frac{1}{2m} \leq N+1$, where

$$p_n(x) = m^n x^n + \cdots + \frac{m!}{(m-n)!} x$$

is a polynomial of degree $n$ without constant term.

We pause to remark that the functions

$$E^{\gamma,\delta}_{\alpha,\beta}(z) := \frac{\Gamma(\delta)}{\Gamma(\gamma)} \sum_{k=0}^{\infty} \frac{\Gamma(k+\gamma) z^k}{\Gamma(\alpha k + \beta) \Gamma(k+\delta)},$$

of which $E^{(n-1)}_{1/m,1/m}$ is the special case $E^{n,1}_{1/m,1/m}$, were recently studied in [19]; however, there seem to be no results on their asymptotic behaviour readily available in the literature.

Recall that $d\sigma$ denotes the surface measure on the unit sphere $\mathbf{S}^{2n-1}$ of $\mathbf{C}^n$.

**Lemma 15.** *As $r|y| \to +\infty$,*

$$\int_{\mathbf{S}^{2n-1}} |K_{\beta m}(r\zeta, y)| \, d\sigma(\zeta) \sim (r|y|)^{\frac{m}{2}-n} e^{\beta r^m |y|^m}.$$

Note that a brute force application of (18) gives only $(m-1)n > \frac{m}{2} - n$ in the exponent instead of $\frac{m}{2} - n$; this improvement is crucial for the applications below.

*Proof.* Consider first the integral

$$I(R) := \int_{-\pi}^{\pi} |E^{(n-1)}_{\frac{1}{m},\frac{1}{m}}(Re^{i\theta})| \, d\theta.$$

For $m > \frac{1}{2}$, the contribution from $\frac{\pi}{2m} < \theta < \pi$ stays bounded as $R \to +\infty$; thus

$$I(R) \sim R^{(m-1)n} \int_{-\pi/2m}^{\pi/2m} e^{R^m \cos m\theta} \, d\theta.$$

For $0 < m \leq \frac{1}{2}$, the term $j = 0$ in (22) always dominates the others, so

$$I(R) \sim R^{(m-1)n} \int_{-\pi}^{\pi} e^{R^m \cos m\theta} \, d\theta.$$

The last two integrals can both be handled by the Laplace method. Recall that the latter asserts that for real-valued functions $f, g$ which are $C^\infty$ on a finite interval $[a, b]$ and such that $g$ attains its maximum at a unique interior point $x_0 \in (a, b)$, $g''(x_0) \neq 0$, the integral

$$\text{(23)} \qquad \int_a^b f(x) \, e^{\lambda g(x)} \, dx$$

behaves for $\lambda \to +\infty$ as

$$\text{(24)} \qquad f(x_0) e^{\lambda g(x_0)} \sqrt{\frac{-2\pi}{\lambda g''(x_0)}}.$$

See e.g. [10], Chapter II, Theorem 1.3.[2] In our case, $\lambda = R^m$, $f \equiv 1$, $g(x) = \cos mx$ and $x_0 = 0$; thus (24) gives that, both for $0 < m \leq \frac{1}{2}$ and $m > \frac{1}{2}$,

$$\text{(25)} \qquad I(R) \sim R^{(m-1)n} e^{R^m} / R^{m/2}.$$

By (17), we have

$$\int_{\mathbf{S}^{2n-1}} |K_{\beta m}(r\zeta, y)| \, d\sigma(\zeta) = \frac{m\beta^{m/n}}{\pi^n} \int_{\mathbf{S}^{2n-1}} |E_{\frac{1}{m},\frac{1}{m}}^{(n-1)}(\beta^{1/m} r\langle \zeta, y\rangle)| \, d\sigma(\zeta).$$

Since the measure $d\sigma$ is invariant under rotations, we may assume without loss of generality that $y = (|y|, 0, \ldots, 0)$; the last integral thus equals

$$\int_{\mathbf{S}^{2n-1}} |E_{\frac{1}{m},\frac{1}{m}}^{(n-1)}(\beta^{1/m} r|y|\zeta_1)| \, d\sigma(\zeta).$$

Parameterizing the sphere as $\text{Im } \zeta_n = \pm\sqrt{1 - |\zeta_1|^2 - |\zeta'|^2}$, where $\zeta' := (\zeta_2, \ldots, \zeta_{n-1}, \text{Re } \zeta_n) \in \mathbf{R}^{2n-3}$, $|\zeta'|^2 < 1 - |\zeta_1|^2$, we have $d\sigma(\zeta) = (1 - |\zeta_1|^2 - |\zeta'|^2)^{-1/2} \, d\zeta_1 \, d\zeta'$. Carrying out the $\zeta'$ integration shows that, up to an immaterial constant factor (involving the volume of the unit ball in $\mathbf{R}^{2n-3}$) which we ignore, the last integral equals

$$\int_{|\zeta_1|<1} |E_{\frac{1}{m},\frac{1}{m}}^{(n-1)}(\beta^{1/m} r|y|\zeta_1)| \, (1 - |\zeta_1|^2)^{n-2} \, d\zeta_1$$

$$= \int_0^1 \int_{-\pi}^{\pi} |E_{\frac{1}{m},\frac{1}{m}}^{(n-1)}(\beta^{1/m} r|y|\rho e^{i\theta})| \, (1 - \rho^2)^{n-2} \, \rho \, d\theta \, d\rho.$$

Substituting (25) for the interior integral yields

$$\int_0^1 (r|y|\rho)^{(m-1)n - \frac{m}{2}} e^{\beta r^m |y|^m \rho^m} (1 - \rho^2)^{n-2} \rho \, d\rho =: \mathcal{Q}.$$

---

[2]Or [9], Section 1.3 (which contains only a sketch of the proof, however).

This is again susceptible to the Laplace method, this time to the variant concerning the integral (23) when $g$ attains its maximum at the boundary point $x_0 = a$ and $g'(a) \neq 0$. In that case, (23) behaves for $\lambda \to +\infty$ as

$$e^{\lambda g(a)} \sum_{k=0}^{\infty} c_k \lambda^{-k-1}$$

where

$$c_k = \left( \frac{1}{-g'(x)} \frac{d}{dx} \right)^k \frac{f(x)}{-g'(x)} \bigg|_{x=a} ;$$

see [10], Chapter II, Theorem 1.1.[3] In our case $[a,b] = [0,1]$, $x = 1-\rho$, $\lambda = \beta r^m |y|^m$, $g(x) = (1-x)^m$, and $f(x) = (r|y|(1-x))^{(m-1)n-\frac{m}{2}}(1-x)(2-x)^{n-2}x^{n-2}$ vanishes at $x = 0$ to order $n-2$. Thus

$$\mathcal{Q} \sim (r|y|)^{(m-1)n-\frac{m}{2}} e^{\beta r^m |y|^m} / (r|y|)^{(n-1)m}.$$

Putting everything together and noting that

$$(m-1)n - \frac{m}{2} - (n-1)m = \frac{m}{2} - n$$

gives the lemma. $\square$

For $\beta > 0$ and $1 \leq p \leq \infty$, define the spaces

$$\mathcal{L}^p_{\beta m} := \{f : e^{-\beta|z|^{2m}/2} f(z) \in L^p(\mathbf{C}^n)\}$$

with the obvious norms. For $p \neq \infty$, we thus have $\mathcal{L}^p_{\beta m} = L^p_{\beta p/2, m}$.

The proof below was inspired by [16].

*Proof of Theorem 1.* We already know that $P_{\beta m} : L^p_{\alpha m} \to L^q_{\gamma m}$ is bounded if $c = c_{p,q}(\alpha, \beta, \gamma) > 1$ (Proposition 9), and unbounded if $c < 1$ (Proposition 7) or $c = 1$ and $p > q$ (Corollary 11; recall that $0 < m \leq 1$ by hypothesis). We thus only need to show that it is bounded if $c = 1$ and $p \leq q$.

By Proposition 6, it is further enough to exhibit only one triple $(\alpha, \beta, \gamma)$ with $c_{p,q}(\alpha, \beta, \gamma) = 1$ for which $P_{\beta m} : L^p_{\alpha m} \to L^q_{\gamma m}$ is bounded. We prove this for

$$\alpha = \frac{\beta p}{2}, \quad \gamma = \frac{\beta q}{2}$$

(for which, indeed, $c = \frac{4\gamma}{\beta^2 q}(\beta - \frac{\alpha}{p}) = \frac{2\beta q}{\beta^2 q} \frac{\beta}{2} = 1$). In other words, we prove that

(26) $$P_{\beta m} : \mathcal{L}^p_{\beta m} \to \mathcal{L}^q_{\beta m} \quad \text{boundedly}$$

for all $1 \leq p \leq q \leq \infty$.

Now by the definition of $\mathcal{L}^p_{\beta m}$, the operator $E_{\beta m}$ from Section 3 maps $L^p(\mathbf{C}^n)$ isometrically onto $\mathcal{L}^p_{\beta m}$, for all $1 \leq p \leq \infty$. Since the spaces $L^p(\mathbf{C}^n)$, $1 \leq p \leq \infty$, form an interpolation scale, so do $\mathcal{L}^p_{\beta m}$. Using interpolation, we therefore conclude that it is enough to establish (26) for the pairs $(p,q) = (1,1), (\infty, \infty)$ and $(1, \infty)$.

---

[3]Or [9], Section 1.2 (again with some details of the proof omitted).

Assume first that $f \in \mathcal{L}^\infty_{\beta m}$. Then

$$|P_{\beta m} f(z)| = \Big| \int_{\mathbf{C}^n} K_{\beta m}(z,w) \, f(w) \, e^{-\beta |w|^{2m}} \, dw \Big|$$
$$\leq \|f\|_{\mathcal{L}^\infty_{\beta m}} \int_{\mathbf{C}^n} |K_{\beta m}(z,w)| \, e^{-\frac{\beta}{2}|w|^{2m}} \, dw.$$

Now by Lemma 15,

$$\int_{\mathbf{C}^n} |K_{\beta m}(z,w)| \, e^{-\frac{\beta}{2}|w|^{2m}} \, dw$$
$$= \int_0^\infty \int_{\mathbf{S}^{2n-1}} |K_{\beta m}(z, r\zeta)| \, e^{-\beta r^{2m}/2} \, r^{2n-1} \, d\sigma(\zeta) \, dr$$
$$\lesssim \int_0^\infty (r|z|)^{\frac{m}{2}-n} \, e^{\beta r^m |z|^m - \frac{\beta}{2} r^{2m}} \, r^{2n-1} \, dr,$$

while by Lemma 8 (or rather its proof)

$$\int_0^\infty r^{\frac{m}{2}-n} \, e^{\beta r^m |z|^m - \frac{\beta}{2} r^{2m}} \, r^{2n-1} \, dr \lesssim |z|^{(\frac{m}{2}-n)+2n-m} \, e^{\beta |z|^{2m}/2};$$

thus

(27) $$\int_{\mathbf{C}^n} |K_{\beta m}(z,w)| \, e^{-\frac{\beta}{2}|w|^{2m}} \, dw \lesssim |z|^{2(\frac{m}{2}-n)+2n-m} e^{\frac{\beta}{2}|z|^{2m}} = e^{\frac{\beta}{2}|z|^{2m}}.$$

Hence
$$|P_{\beta m} f(z)| \lesssim \|f\|_{\mathcal{L}^\infty_{\beta m}} \, e^{\frac{\beta}{2}|z|^{2m}},$$

proving that $P_{\beta m}$ maps $\mathcal{L}^\infty_{\beta m}$ boundedly into itself.

Assume next that $f \in \mathcal{L}^1_{\beta m}$. Then

$$\|P_{\beta m} f\|_{\mathcal{L}^1_{\beta m}} = \int_{\mathbf{C}^n} \Big| \int_{\mathbf{C}^n} K_{\beta m}(z,w) \, f(w) \, e^{-\beta|w|^{2m}} \, dw \Big| e^{-\frac{\beta}{2}|z|^{2m}} \, dz$$
$$\leq \int_{\mathbf{C}^n} |f(w)| \int_{\mathbf{C}^n} |K_{\beta m}(z,w)| \, e^{-\frac{\beta}{2}|z|^{2m}} \, dz \, e^{-\beta|w|^{2m}} \, dw$$
$$\lesssim \int_{\mathbf{C}^n} |f(w)| \, e^{-\frac{\beta}{2}|w|^{2m}} \, dw = \|f\|_{\mathcal{L}^1_{\beta m}}$$

by (27). Thus $P_{\beta m}$ maps $\mathcal{L}^1_{\beta m}$ boundedly into itself.

Finally, note that for $0 < m \leq 1$, it follows from (18) that

$$|K_{\beta m}(z,w)| \lesssim e^{\beta |z|^m |w|^m} \leq e^{\beta \frac{|z|^{2m}+|w|^{2m}}{2}}.$$

Thus for any $f \in \mathcal{L}^1_{\beta m}$

$$|P_{\beta m} f(z)| \, e^{-\frac{\beta}{2}|z|^{2m}} = \Big| \int_{\mathbf{C}^n} K_{\beta m}(z,w) \, f(w) \, e^{-\beta|w|^{2m} - \frac{\beta}{2}|z|^{2m}} \, dw \Big|$$
$$\lesssim \int_{\mathbf{C}^n} |f(w)| \, e^{-\frac{\beta}{2}|w|^{2m}} \, dw = \|f\|_{\mathcal{L}^1_{\beta m}},$$

so $P_{\beta m}$ maps $\mathcal{L}^1_{\beta m}$ boundedly into $\mathcal{L}^\infty_{\beta m}$, and the proof of Theorem 1 is complete. □

Note that the only place where the hypothesis $0 < m \le 1$ was used was the very last part of the proof — when showing that $P_{\beta m}$ is bounded from $\mathcal{L}^1_{\beta m}$ into $\mathcal{L}^\infty_{\beta m}$. In particular, (26) holds for $(p,q) = (1,1)$ and $(\infty,\infty)$ for all $m > 0$; by interpolation, we thus have (26) for all $1 \le p = q < \infty$. By Proposition 6, we thus obtain the following result.

**Proposition 16.** $P_{\beta m} : L^p_{\alpha m} \to L^q_{\gamma m}$ is bounded whenever $c = 1$ and $p = q$, i.e. when $p = q = q_{\max}$.

We are now ready to prove the last item that remains to be verified.

*Proof of Theorem 2.* Since we already know $P_{\beta m} : L^p_{\alpha m} \to L^q_{\gamma m}$ to be bounded when $c > 1$ (Proposition 9) and unbounded when $c < 1$ (Proposition 7) or $c = 1$ and $p \ne q$ (Corollary 14), the claim is immediate from the last proposition. □

## 8. Concluding remarks

**8.1.** Using the method behind Lemma 15, it is possible to give another proof of the necessity criterion from Section 4 (Proposition 7), avoiding the use of Gamma functions and Stirling's formula. Namely, let $\chi_a$, $a \in \mathbf{C}^n$, be the characteristic function of the polydisc of radius 1 with center $a$, normalized to be of total mass one (i.e. multiplied by $\pi^{-n}$). From the mean value property of holomorphic functions

$$f(a) = \int_{\mathbf{C}^n} f \chi_a \, dz = \int_{\mathbf{C}^n} f(z) \chi_a(z) e^{\beta |z|^{2m}} \, d\mu_{\beta m}(z)$$

and the reproducing property of $K_{\beta m}$ it follows that $K_{\beta m}(\,\cdot\,, a) = P_{\beta m}(\chi_a e^{\beta |\cdot|^{2m}})$. Thus if $P_{\beta m} : L^p_{\alpha m} \to L^q_{\gamma m}$ is bounded, then

$$\|K_{\beta m}(\,\cdot\,, a)\|_{q\gamma m} \le C \,\|\chi_a e^{\beta |\cdot|^{2m}}\|_{p\alpha m} \qquad \forall a \in \mathbf{C}^n$$

with some $C$ independent of $a$. Proceeding as in the proof of Lemma 15 and using Lemma 8 (or rather its proof) yields

$$\|K_{\beta m}(\,\cdot\,, a)\|_{q\gamma m} \sim |a|^{2n(m-1)(1-\frac{1}{q})} e^{\beta^2 q |a|^{2m}/4\gamma}.$$

On the other hand

$$\|\chi_a e^{\beta |\cdot|^{2m}}\|_{p\alpha m} = \left[\int_{\mathbf{C}^n} \chi_a(z) \, e^{(\beta p - \alpha)|z|^{2m}} \, dz\right]^{\frac{1}{p}} \lesssim e^{(\beta - \frac{\alpha}{p})(|a|+1)^{2m}}.$$

Letting $a \to \infty$ we thus see that $\frac{\beta^2 q}{4\gamma} \le \beta - \frac{\alpha}{p}$, or $1 \le c$.

Unfortunately, this method — apart from its dependence on the fairly sophisticated machinery needed in the proof of Lemma 15 — seems to have the drawback of being too coarse to reproduce also Propositions 10 and 12 in the borderline case $c = 1$.

**8.2.** Another proof of Proposition 16 can be given using the Schur test (pretty much like in §3 of [6]); however, we do not know how to use the Schur test to prove also Theorem 1. Here are the details.

In its strengthened variant, taken from Gagliardo [12] (where it is attributed to Aronszajn [2]), the Schur test runs as follows. Let $(X, \mu)$ and $(Y, \nu)$ be two totally $\sigma$-finite measure spaces, $K(x, y) \geq 0$ a nonnegative measurable function on $X \times Y$, and consider the linear transformations

$$(Tu)(y) = \int_X K(x,y) u(x) \, d\mu(x), \qquad (T^*v)(x) = \int_Y K(x,y) v(y) \, d\nu(y).$$

Let $1 < q \leq p < +\infty$, $\frac{1}{p'} + \frac{1}{p} = \frac{1}{q'} + \frac{1}{q} = 1$. Then a sufficient condition in order that $T$ restricted to $L^p(X, \mu)$ be a bounded operator from $L^p(X, \mu)$ into $L^q(Y, \nu)$ is that there exist two measurable functions $\phi(x)$, $\psi(y)$, positive and finite a.e., and a finite constant $C$ such that

(28) $$T\phi \leq C \, \psi^{q'/q},$$

(29) $$T^*\psi \leq C \, \phi^{p/p'},$$

and

(30) $$\iint_{X \times Y} K(x,y) \phi(x) \psi(y) \, d\mu(x) \, d\nu(y) \leq C \qquad \text{if } p \neq q.$$

For the proof see [12], pp. 429–430.

We apply the test to $X = Y = \mathbf{C}^n$, $d\mu = d\mu_{\alpha m}$, $d\nu = d\mu_{\gamma m}$, and $K(x,y) = |K_{\beta m}(x,y)| e^{(\alpha-\beta)|x|^{2m}}$, with the choice

$$\phi(x) = e^{\lambda p' |x|^{2m}}, \quad \psi(y) = e^{\nu q |y|^{2m}},$$

where $\lambda$ and $\nu$ are to be specified later on. Using Lemma 15 and Lemma 8 (or rather its proof), one can show (cf. the derivation of (27) above) that for $C > 0$,

$$\int_{\mathbf{C}^n} |K_{\beta m}(x,y)| \, e^{-C|x|^{2m}} \, dx \sim \int_0^\infty e^{-Cr^{2m}} (r|y|)^{\frac{m}{2}-n} e^{\beta r^m |y|^m} r^{2n-1} \, dr$$

(31) $$\sim e^{\beta^2 |y|^{2m}/4C}.$$

Using this, the conditions (28) and (29) are seen to be equivalent, respectively, to

(32) $$\beta - \lambda p' > 0 \quad \text{and} \quad \frac{\beta^2}{4(\beta - \lambda p')} \leq \nu q',$$

(33) $$\gamma - \nu q > 0 \quad \text{and} \quad \alpha - \beta + \frac{\beta^2}{4(\gamma - \nu q)} \leq \lambda p,$$

while (30) becomes simply

(34) $$\beta^2 < 4(\beta - \lambda p')(\gamma - \nu q).$$

For $p = q = q_{\max}$ (so that we need not worry about (30)), the conditions (32) and (33) are satisfied for $\nu = \gamma/qq'$ and $\lambda = \alpha/pp'$. This gives Proposition 16.

Additionally, we can get from the above also a proof of Proposition 9 for $p \geq q$ (albeit much more complicated than the one in Section 5). Indeed, the inequalities (32) and (33) can be equivalently rewritten, respectively, as

$$(35) \qquad 0 < \frac{\beta^2}{4\nu q' p'} \leq \frac{\beta}{p'} - \lambda, \qquad 0 < \frac{\beta^2}{4(\gamma - \nu q)p} \leq \frac{\beta - \alpha}{p} + \lambda.$$

The left halves of these inequalities mean simply that $0 < \nu < \frac{\gamma}{q}$, so we set

$$(36) \qquad \nu = (1-x)\frac{\gamma}{q}, \quad 0 < x < 1.$$

Their right halves are then solvable for $\lambda$ if and only if they fulfill the compatibility condition

$$\frac{\beta^2}{4\nu q' p'} + \frac{\beta^2}{4(\gamma - \nu q)p} \leq \frac{\beta}{p'} + \frac{\beta - \alpha}{p} = \beta - \frac{\alpha}{p}.$$

Using (36) this becomes, after a small manipulation,

$$\frac{1}{(1-x)p'q'} + \frac{1}{xpq} \leq c,$$

or

$$(37) \qquad \left(x - \frac{1}{p}\right)\left(x - \frac{1}{q}\right) \leq (c-1)x(1-x).$$

For $p = q$, we thus see that $x = \frac{1}{p} = \frac{1}{q}$ (which, by the way, corresponds to the choice $\nu = \lambda = 1/pp' = 1/qq'$ we used in the preceding paragraph) does the job (for any $c > 1$). For $p > q$, the extra condition (34) can be rewritten as

$$0 < \frac{\beta^2}{4(\gamma - \nu q)p'} < \frac{\beta}{p'} - \lambda,$$

which in conjunction with (35) is equivalent to one more compatibility condition

$$\frac{\beta^2}{4(\gamma - \nu q)p'} + \frac{\beta^2}{4(\gamma - \nu q)p} < \frac{\beta}{p'} + \frac{\beta - \alpha}{p} = \beta - \frac{\alpha}{p} = \frac{\beta^2 qc}{4\gamma},$$

or, upon a small manipulation,

$$(38) \qquad \frac{1}{x} < qc, \quad \text{or} \quad \frac{1}{qc} < x.$$

Thus for $c > 1$, $x = 1/q$ again solves both (37) and (38), hence yielding a solution $\nu = 1/qq'$ and $\lambda$ to (32), (33) and (34). Thus we have recovered Proposition 7 for $p \leq q$.

Unfortunately, this approach no longer works when $p > q$ and $c = 1$. Indeed, in that case (37) is solved by $x \in [\frac{1}{p}, \frac{1}{q}]$, while (38) by $x \in (\frac{1}{q}, 1)$, so there are no common solutions.

Taking more general test functions $\phi, \psi$ of the form
$$|x|^A e^{B|x|^m - C|x|^{2m}}$$
(and using the appropriate generalization
$$\int_{\mathbf{C}^n} |K_{\beta m}(x,y)| \, |x|^A \, e^{B|x|^m - C|x|^{2m}} \, dx \lesssim |y|^A \, e^{\frac{2\beta B|y|^m + \beta^2 |y|^{2m}}{4C}}$$
of (31)) does not help.

Finally, there is no variant of the Schur test that would be applicable for $p < q$ (see the discussion in [12] for clarification). In particular, other methods need probably to be used for the proofs of e.g. Proposition 9 when $p < q$.

**8.3.** Some of our proofs — e.g. those in Sections 2, 3 and 6 — did not really require that $p, q$ be greater or equal to 1. The corresponding results thus remain in force for all $p, q > 0$.

**8.4.** Similarly, we have made no real effort to analyze also the case of $p, q = \infty$. It is quite conceivable that most of the results extend also to this case, possibly with $L^\infty_{\alpha m} = L^\infty(\mathbf{C}^n)$ replaced by the spaces $\mathcal{L}^\infty_{\alpha m}$ from Section 7.

**8.5.** A possible way of producing some further necessary conditions (which might be needed for the remaining cases of $1 < m \leq \frac{2n}{2n-1}$ and $m \geq 2$), like those in Sections 4 and 6, is testing the boundedness of $P_{\beta m}$ on more complicated test functions than $f(z) = z^\nu e^{\lambda |z|^{2m}}$. A good candidate might be
$$f(z) = \overline{z}^\rho z^{\nu + \rho} \phi(|z|),$$
with the analogue of (8) given by
$$P_{\beta m} f(w) = \frac{w^\nu}{\|z^\nu\|^2_{\beta m}} \int_{\mathbf{C}^n} |z^{\rho+\nu}|^2 \, \phi(|z|) \, d\mu_{\beta m}(z)$$
$$= \frac{w^\nu}{\|z^\nu\|^2_{\beta m}} \frac{2\pi^n (\rho + \nu)!}{\Gamma(n + |\rho| + |\nu|)} \int_0^\infty r^{2|\rho+\nu|} \phi(r) \, e^{-\beta r^{2m}} \, r^{2n-1} \, dr.$$
For $\phi$ of the form $\phi(r) = r^A e^{Br^m + Cr^{2m}}$, the last integral can be evaluated explicitly, or at least estimated asymptotically as $\nu$ (or perhaps other parameters?) gets large. (See e.g. Evgrafov [8], §1.2, especially Theorems 1.2.5 and 1.2.6.)

**8.6.** At first sight, there is also an avenue left open in our proof of Proposition 9 which could yield additional sufficient conditions for the borderline case $c = 1$. Namely, the integral in (19) there can still be finite even if $c = 1$ (so that the two exponential terms cancel out), but the power of $|z|$ decays sufficiently fast:
$$[2(m-1)np' + 2n - m]\frac{q}{p'} < -2n.$$
Upon a small manipulation the last inequality assumes the form
$$m < \frac{2n(\frac{1}{p} - \frac{1}{q})}{2n - \frac{1}{p'}}.$$
However, since
$$\frac{2n(\frac{1}{p} - \frac{1}{q})}{2n - \frac{1}{p'}} = \frac{2n(\frac{1}{p} - \frac{1}{q})}{2n - 1 + \frac{1}{p}} \leq \frac{2n(\frac{1}{p} - \frac{1}{q})}{\frac{2n-1}{p} + \frac{1}{p}} = \frac{\frac{1}{p} - \frac{1}{q}}{\frac{1}{p}} < 1,$$
this can only apply in the case of $0 < m \leq 1$ and $p < q$, which is already fully covered by Theorem 1. Thus we do not get anything new.

**8.7.** We have seen in Theorems 1 and 2 that for $0 < m \leq 1$ and $\frac{2n}{2n-1} < m < 2$, the boundedness of $P_\beta : L^p_{\alpha m} \to L^q_{\gamma m}$ depends only on the signs of $c - 1 = \frac{q_{\max}}{q} - 1$ and $p - q$. It can be shown that this is in fact the case in general.

**Proposition 17.** *For each fixed $m$ and $n$, there are only four (mutually exclusive) possibilities that can occur:*

  (a) $P_{\beta m} : L^p_{\alpha m} \to L^q_{\gamma m}$ *is bounded iff either $c > 1$, or $c = 1$ and $p = q$;*
  (b) $P_{\beta m} : L^p_{\alpha m} \to L^q_{\gamma m}$ *is bounded iff either $c > 1$, or $c = 1$ and $p \leq q$;*
  (c) $P_{\beta m} : L^p_{\alpha m} \to L^q_{\gamma m}$ *is bounded iff either $c > 1$, or $c = 1$ and $p \geq q$;*
  (d) $P_{\beta m} : L^p_{\alpha m} \to L^q_{\gamma m}$ *is bounded iff $c \geq 1$.*

*Proof.* Since $P^2_{\beta m} = P_{\beta m}$ ($P_{\beta m}$ being a projection), the boundedness of $P_{\beta m} : L^p_{\alpha m} \to L^q_{\gamma m}$ will follow if we prove that

(39) $$P_{\beta m} : L^p_{\alpha m} \to L^r_{\delta m}, \quad L^r_{\delta m} \to L^q_{\gamma m} \qquad \text{(boundedly)}$$

for some $r \geq 1$, $\delta > 0$. Furthermore, we are interested only in the case when $c := \frac{4\gamma}{\beta^2 q}(\beta - \frac{\alpha}{p}) = 1$, since we already know that $P_{\beta m} : L^p_{\alpha m} \to L^q_{\gamma m}$ is bounded if $c > 1$ and unbounded if $c < 1$.

As in the proof of Theorem 1 in Section 7, owing to Proposition 6 it is enough to deal with the special values

$$\alpha = \frac{\beta p}{2}, \quad \gamma = \frac{\beta q}{2}.$$

Since $L^p_{\beta p/2, m} = \mathcal{L}^p_{\beta m}$, (39) becomes simply

$$P_{\beta m} : \mathcal{L}^p_{\beta m} \to \mathcal{L}^r_{\beta m}, \quad \mathcal{L}^r_{\beta m} \to \mathcal{L}^q_{\beta m} \qquad \text{(boundedly)}.$$

We fix $m, n$ and $\beta$ from now on. Let us denote

$$A^\infty := \{(\tfrac{1}{p}, \tfrac{1}{q}) \in [0,1] \times [0,1] : P_{\beta m} : \mathcal{L}^p_{\beta m} \to \mathcal{L}^q_{\beta m} \text{ boundedly}\},$$
$$A := A^\infty \cap (0,1] \times (0,1].$$

Then by interpolation (since $\mathcal{L}^p_{\beta m}$ form an interpolation scale)

(40) $$A^\infty \text{ is convex,}$$

and (by Proposition 16) $A^\infty$ contains the diagonal:

(41) $$\{(x,x) : 0 \leq x \leq 1\} \subset A^\infty.$$

On the other hand, (39) translates into

(42) $$(x,y), (y,z) \in A^\infty \implies (x,z) \in A^\infty.$$

Assume that $A \ni (x,y)$ with some $x > y$. By (40) and (41), the segment connecting $(x,y)$ and $(0,0)$ lies in $A$, i.e.

$$\{(tx, ty) : 0 < t \leq 1\} \subset A.$$

Taking $t = y/x$, (42) gives $(x, y^2/x) \in A$. Taking then $t = y^2/x^2$, (42) gives $(x, y^3/x^2) \in A$. Continuing, we get $(x, (y/x)^k y) \in A$ for all $k = 1, 2, 3, \ldots$, and by (40) we see that

$$(x, z) \in A \quad \forall\, 0 < z \leq x.$$

Finally, also the segment connecting $(x, y)$ with $(1, 1)$ lies in $A$ by (40), and we can apply the above reasoning to any point of this segment in the place of $(x, y)$, leading to the conclusion that

$$(w, z) \in A \quad \forall\, 0 < z \leq w \leq 1.$$

In a completely similar fashion one can show that if $A \ni (x, y)$ with some $x < y$, then

$$(w, z) \in A \quad \forall\, 0 < w \leq z \leq 1.$$

Summarizing, we see that there are only four possible cases:
- $A = \{(x, y) \in (0, 1] \times (0, 1] : x = y\}$,
- $A = \{(x, y) \in (0, 1] \times (0, 1] : x \geq y\}$,
- $A = \{(x, y) \in (0, 1] \times (0, 1] : x \leq y\}$,
- $A = (0, 1] \times (0, 1]$.

Translating this back into terms of $p, q, \alpha, \beta, \gamma$, the claim in the proposition follows. $\square$

From our results we know that the possibility (d) never occurs, (b) occurs for $0 < m \leq 1$, and (a) for $\frac{2n}{2n-1} < m < 2$. The "wild guesses" mentioned after Corollary 14 are that (a) prevails in fact for all $1 < m < 2$, and (c) for $m \geq 2$. Although we do not know how to prove this, it transpires from the above that it might be enough to understand the extreme cases $P_{\beta m} : \mathcal{L}^1_{\beta m} \to \mathcal{L}^\infty_{\beta m}, \mathcal{L}^\infty_{\beta m} \to \mathcal{L}^1_{\beta m}$.

(H.B. and H.Y.) LATP, U.M.R. C.N.R.S. 6632, CMI, Université de Provence, 39 Rue F-Joliot-Curie, 13453 Marseille Cedex 13, France
*E-mail address*: `bommier@gyptis.univ-mrs.fr, youssfi@gyptis.univ-mrs.fr`

(M.E.) Mathematics Institute, Žitná 25, 11567 Prague 1, Czech Republic and Mathematics Institute, Silesian University at Opava, Na Rybníčku 1, 74601 Opava, Czech Republic
*E-mail address*: `englis@math.cas.cz`